\catcode`\^^Z=9
\catcode`\^^M=10
\output={\if N\header\headline={\hfill}\fi
\plainoutput\global\let\header=Y}
\magnification\magstep1
\tolerance = 500
\hsize=14.4true cm
\vsize=22.5true cm
\parindent=6true mm\overfullrule=2pt
\newcount\kapnum \kapnum=0
\newcount\parnum \parnum=0
\newcount\procnum \procnum=0
\newcount\nicknum \nicknum=1
\font\ninett=cmtt9

\font\ninebf=cmbx9

\font\sixbf=cmbx6
\font\ninesl=cmsl9

\font\nineit=cmti9

\font\ninerm=cmr9

\font\sixrm=cmr6
\font\ninei=cmmi9
\font\eighti=cmmi8
\font\sixi=cmmi6
\skewchar\ninei='177 \skewchar\eighti='177 \skewchar\sixi='177
\font\ninesy=cmsy9
\font\eightsy=cmsy8
\font\sixsy=cmsy6
\skewchar\ninesy='60 \skewchar\eightsy='60 \skewchar\sixsy='60
\font\titelfont=cmr10 scaled 1440
\font\paragratit=cmbx10 scaled 1200

\font\name=cmcsc10
\font\emph=cmbxti10

\font\tenmsbm=msbm10
\font\sevenmsbm=msbm7
%

%

%
\font\teneufm=eufm10
\font\seveneufm=eufm7
\font\fiveeufm=eufm5
\newfam\eufmfam
\textfont\eufmfam=\teneufm
\scriptfont\eufmfam=\seveneufm
\scriptscriptfont\eufmfam=\fiveeufm

\font\tenmsam=msam10
\font\sevenmsam=msam7
\font\fivemsam=msam5
\newfam\msamfam
\textfont\msamfam=\tenmsam
\scriptfont\msamfam=\sevenmsam
\scriptscriptfont\msamfam=\fivemsam
\font\tenmsbm=msbm10
\font\sevenmsbm=msbm7
\font\fivemsbm=msbm5
\newfam\msbmfam
\textfont\msbmfam=\tenmsbm
\scriptfont\msbmfam=\sevenmsbm
\scriptscriptfont\msbmfam=\fivemsbm
\def\Bbb#1{{\fam\msbmfam\relax#1}}
\def\cz{{\kern0.4pt\Bbb C\kern0.7pt}
}
\def\ez{{\kern0.4pt\Bbb E\kern0.7pt}
}
\def\fz{{\kern0.4pt\Bbb F\kern0.3pt}}
\def\gz{{\kern0.4pt\Bbb Z\kern0.7pt}}
\def\hz{{\kern0.4pt\Bbb H\kern0.7pt}
}
\def\kz{{\kern0.4pt\Bbb K\kern0.7pt}
}
\def\nz{{\kern0.4pt\Bbb N\kern0.7pt}
}
\def\oz{{\kern0.4pt\Bbb O\kern0.7pt}
}
\def\rz{{\kern0.4pt\Bbb R\kern0.7pt}
}
\def\sz{{\kern0.4pt\Bbb S\kern0.7pt}
}
\def\pz{{\kern0.4pt\Bbb P\kern0.7pt}
}
\def\qz{{\kern0.4pt\Bbb Q\kern0.7pt}
}
\newskip\ttglue
\def\ninepoint{\def\rm{\fam0\ninerm}%
  \textfont0=\ninerm \scriptfont0=\sixrm \scriptscriptfont0=\fiverm
  \textfont1=\ninei \scriptfont1=\sixi \scriptscriptfont1=\fivei
  \textfont2=\ninesy \scriptfont2=\sixsy \scriptscriptfont2=\fivesy
  \textfont3=\tenex \scriptfont3=\tenex \scriptscriptfont3=\tenex
  \def\it{\fam\itfam\nineit}%
  \textfont\itfam=\nineit
  \def\sl{\fam\slfam\ninesl}%
  \textfont\slfam=\ninesl
  \def\bf{\fam\bffam\ninebf}%
  \textfont\bffam=\ninebf \scriptfont\bffam=\sixbf
   \scriptscriptfont\bffam=\fivebf
  \def\tt{\fam\ttfam\ninett}%
  \textfont\ttfam=\ninett
  \tt \ttglue=.5em plus.25em minus.15em
  \normalbaselineskip=11pt
  \font\name=cmcsc9
  \let\sc=\sevenrm
  \let\big=\ninebig
  \setbox\strutbox=\hbox{\vrule height8pt depth3pt width0pt}%
  \normalbaselines\rm
  \def\sl{\it}}

\headline={\ifodd\pageno\rightheadline\else\leftheadline\fi}
\def\rightheadline{\ninepoint Paragraphen"uberschrift\hfill\folio}
\def\leftheadline{\ninepoint\folio\hfill Chapter"uberschrift}
\let\header=Y
\def\titel#1{\need 9cm \vskip 2truecm
\parnum=0\global\advance \kapnum by 1
{\baselineskip=16pt\lineskip=16pt\rightskip0pt
plus4em\spaceskip.3333em\xspaceskip.5em\pretolerance=10000\noindent
\titelfont Chapter \uppercase\expandafter{\romannumeral\kapnum}.
#1\vskip2true cm}\def\leftheadline{\ninepoint
\folio\hfill Chapter \uppercase\expandafter{\romannumeral\kapnum}.
#1}\let\header=N
}
\def\Titel#1{\need 9cm \vskip 2truecm
\global\advance \kapnum by 1
{\baselineskip=16pt\lineskip=16pt\rightskip0pt
plus4em\spaceskip.3333em\xspaceskip.5em\pretolerance=10000\noindent
\titelfont\uppercase\expandafter{\romannumeral\kapnum}.
#1\vskip2true cm}\def\leftheadline{\ninepoint
\folio\hfill\uppercase\expandafter{\romannumeral\kapnum}.
#1}\let\header=N
}
\def\need#1cm {\par\dimen0=\pagetotal\ifdim\dimen0<\vsize
\global\advance\dimen0by#1 true cm
\ifdim\dimen0>\vsize\vfil\eject\noindent\fi\fi}
\def\neupara#1{\par\penalty-2000
\procnum=0\global\advance\parnum by 1
\vskip1cm\noindent{\paragratit \the\parnum. #1}%
\def\rightheadline{\ninepoint\S\the\parnum.\ #1\hfill \folio}%
\vskip 8mm\noindent}
\def\Proclaim #1 #2\finishproclaim {\bigbreak\noindent
{\bf#1\unskip{}. }{\it#2}\medbreak\noindent}
%
\gdef\proclaim #1 #2 #3\finishproclaim {\bigbreak\noindent%
\global\advance\procnum by 1
{%
{\relax\ifodd \nicknum
\hbox to 0pt{\vrule depth 0pt height0pt width\hsize
   \quad \ninett#3\hss}\else {}\fi}%
\bf\the\parnum.\the\procnum\ #1\unskip{}. }
{\it#2}
\immediate\write\num{\string\def
 \expandafter\string\csname#3\endcsname
 {\the\parnum.\the\procnum}}
\medbreak\noindent}
\newcount\stunde \newcount\minute \newcount\hilfsvar
\def\uhrzeit{
    \stunde=\the\time \divide \stunde by 60
    \minute=\the\time
    \hilfsvar=\stunde \multiply \hilfsvar by 60
    \advance \minute by -\hilfsvar
    \ifnum\the\stunde<10
    \ifnum\the\minute<10
    0\the\stunde:0\the\minute~Uhr
    \else
    0\the\stunde:\the\minute~Uhr
    \fi
    \else
    \ifnum\the\minute<10
    \the\stunde:0\the\minute~Uhr
    \else
    \the\stunde:\the\minute~Uhr
    \fi
    \fi
    }
 \def\calB{{\cal B}}

\def\calE{{\cal E}}

 \def\calX{{\cal X}}

\def\dim{\mathop{\rm dim}\nolimits}

\def\GL{\mathop{\rm GL}\nolimits}

\def\kernel{\mathop{\rm kernel}\nolimits}

\def\U{{\rm U}}

\def\proj{\mathop{\rm proj}\nolimits}

\def\boxit#1{
  \vbox{\hrule\hbox{\vrule\kern6pt
  \vbox{\kern8pt#1\kern8pt}\kern6pt\vrule}\hrule}}
\def\Boxit#1{
  \vbox{\hrule\hbox{\vrule\kern2pt
  \vbox{\kern2pt#1\kern2pt}\kern2pt\vrule}\hrule}}
\def\rahmen#1{
  $$\boxit{\vbox{\hbox{$ \displaystyle #1 $}}}$$}

\def\smallni{\smallskip\noindent }
\def\medni{\medskip\noindent }

\def\lo{\longrightarrow}

\def\loma{\longmapsto}

\def\spitz#1{\langle#1\rangle}
\def\square{\hbox{\hbox to 0pt{$\sqcup$\hss}\hbox{$\sqcap$}}}
\def\qed{\ifmmode\square\else{\unskip\nobreak\hfil
\penalty50\hskip3em\null\nobreak\hfil\square
\parfillskip=0pt\finalhyphendemerits=0\endgraf}\fi}
\def\pn{\the\parnum.\the\procnum}
\def\downmapsto{{\buildrel
        {\vbox{\hbox{\hskip.2pt$\scriptstyle-$}}}
        \over{\raise7pt\vbox{\vskip-4pt\hbox{$\textstyle\downarrow$}}}}}
\openin5=ballcy2.num
\ifeof5 \relax\else
\input ballcy2.num \fi
\nopagenumbers
\immediate\newwrite\num

\let\header=N
\def\transpose#1{\kern1pt{^t\kern-1pt#1}}%

\def\bull{{\hbox{\raise.5ex\hbox{\titelfont.}}}}
\def\bullet{\bull}
 \immediate\openout\num=ballcy2.num
 \immediate\newwrite\num\immediate\openout\num=ballcy2.num
\def\RAND#1{\vskip0pt\hbox to 0mm{\hss\vtop to 0pt{%
  \raggedright\ninepoint\parindent=0pt%
  \baselineskip=1pt\hsize=2cm #1\vss}}\noindent}
\noindent
\centerline{\titelfont Some ball quotients with a Calabi--Yau model}%
\def\leftheadline{\ninepoint\folio\hfill
Some ball quotients with a Calabi--Yau model}%
\def\rightheadline{\ninepoint Introduction\hfill \folio}%
\headline={\ifodd\pageno\rightheadline\else\leftheadline\fi}

\vskip 1.5cm
\leftline{\it \hbox to 6cm{Eberhard Freitag\hss}
Riccardo Salvati
Manni  }
  \leftline {\it  \hbox to 6cm{Mathematisches Institut\hss}
Dipartimento di Matematica, }
\leftline {\it  \hbox to 6cm{Im Neuenheimer Feld 288\hss}
Piazzale Aldo Moro, 2}
\leftline {\it  \hbox to 6cm{D69120 Heidelberg\hss}
 I-00185 Roma, Italy. }
\leftline {\tt \hbox to 6cm{freitag@mathi.uni-heidelberg.de\hss}
salvati@mat.uniroma1.it}
\vskip1cm
\centerline{\paragratit \rm  2012}%
\vskip5mm\noindent%
\let\header=N%
\medni
\neupara{Introduction}%
In a sequence of papers, cf.\ [FS1], [FS2], [FS3] and  [CFS], we studied Siegel 
threefolds which admit a Calabi--Yau model. We got a long list of projective  Calabi--Yau manifolds.
Some of them have been new. 
In this paper we start to look to other kinds of  modular examples.  A natural class are
the ball quotients which belong to the unitary group $\U(1,3)$. Its arithmetic subgroups are
called Picard modular groups and the corresponding varieties are called Picard modular varieties.
In [FS4] we determined explicitly  a  Picard modular variety of  general type. 
On the regular locus of this variety there are many holomorphic three forms with known zero divisors
which have been constructed as Borcherds products.
Resolutions of quotients of this variety, such that the
zero divisors are in the branch locus, are candidates for Calabi--Yau manifolds.
The aim of this note is to treat one distinguished example for this.
In fact  we shall recover a known variety ([Me], Chap.~5, Sect.~6) given by the equations
\rahmen{X_0X_1X_2=X_3X_4X_5,\quad X_0^3+X_1^3+X_2^3=X_3^3+X_4^3+X_5^3.}
as a Picard modular variety with respect to a certain Picard modular group
$G'$. As it has been explained in [Me],
this variety has  a projective small resolution which is a rigid
Calabi-Yau manifold ($h^{12}=0$) with Euler number $72$.  
Moreover a certain group acts on this variety,
fixing the holomorphic differential 3-form. We will show that this group is modular in the sense
that it is given by a subgroup of the normalizer of $G'$ in the full modular group.
\neupara{Differential forms}%
Let $V$ be a finite
dimensional complex vector space which is equipped with a
hermitean form $\spitz{\cdot,\cdot}$ of signature $(1,n)$.
A line (=one dimensional sub-vector-space of $V$) is called positive
if it is  represented by an element of positive norm ($\spitz{z,z}>0$).
We denote by $\calB=\calB(V)\subset P(V)$ the set of all positive
lines. (As usual $P(V)$ denotes the projective space of $V$, i.e.\
the set of all lines in $V$.)
The unitary group $\U(V)\cong\U(1,n)$ acts on $\calB$ as group
of biholomorphic automorphisms.
\smallskip
We choose a sub-vector space $W\in V$ of codimension one which contains no vector of positive
norm. We also consider a vector $e\in V$ such that
$$V=\cz e\oplus W.$$
If $v=Ce+w$, $w\in W$, is a vector of positive norm, then $C\ne 0$. Hence we can normalize
the line $\cz v$ such that $C=1$. This gives an embedding of $\calB$ into $W$.
The linear action of $\U(V)$ on $\calB(V)$ gives a rational action
$z\mapsto g\spitz z$ on the image of $\calB$ in $W$. If we write $g(e+z)=Ce+w$, then
$g\spitz z=w/C$. We have
$$g(e+z)=j(g,z)(e+g\spitz z)\quad\hbox{where}\quad j(g,z):=C.$$
The function $j(g,z)$ is an automorphy factor. We call it the canonical
automorphy factor. Another automorphy factor is the Jacobian.
\proclaim
{Lemma}
{The Jacobian of the transformation
$z\mapsto g\spitz z$ for $g\in \U(V)$ is
$$\det(g)j(g,z)^{-(n+1)}.$$}
CanAut%
\finishproclaim
{\it Proof.\/} The proof can be given for reflections by a somewhat tedious
computation.\qed
\smallskip
Let $M$ be a lattice in $V$, i.e.\ a discrete subgroup with compact quotient. We assume that
$M$ admits complex
multiplication. Then the set of all complex numbers $a\in\cz$ such that $aM\subset M$ is an order
in an imaginary number field. The Picard modular group $\U(M)$ is the subgroup of $\U(V)$ that
preserves $M$.
\smallskip
We recall the notion of a modular form of integral weight $k$. For this we consider the inverse image
$\tilde\calB$ of $\calB$ in $V-\{0\}$. It consists of all $z\in V$ with
$\spitz{z,z}>0$. This is a connected open subset.
\smallskip
Let $G\subset\U(M)$ be a subgroup of finite index and let $v:G\to\cz^\bullet$
be a character.
\eject\noindent
A modular form of weight $k$ on $G$ with respect
to a $v$ is a holomorphic function
$f:\tilde\calB\to\cz$ with the following properties:
\vskip1mm
\item{1)} $f(tz)=t^{-k}f(z)$ for all $t\in\cz^\bullet$,
\item{2)} $f(\gamma z)=v(\gamma)f(z)$ for all $\gamma\in G$,
\item{3)} $f$ is regular at the cusps.
\smallni
We denote this space of modular forms by $[G,k,v]$ or simply by $[G,k]$ when
$v$ is trivial. The ring of modular forms
$$A(G)=\bigoplus[G,k]$$
is a finitely generated algebra. From the theory of Baily and Borel follows that
$A(G)$ is a finitely generated algebra whose associated projective variety can be
identified with the compactification $X(G)$ of $\calB/G$ by (finitely many) cusps,
$$\proj A(G)\cong X(G).$$
From Lemma \CanAut\ we obtain the following result.
\proclaim
{Lemma}
{There is a one-to-one correspondence between $G$-invariant differential forms
of top degree $n$
on $\calB$ and modular forms of weight $n+1$ and multiplier system
$v(\gamma)=\det(\gamma)^{-1}$.}
InvDif%
\finishproclaim
We mention a well-known general fact.
A  $G$-invariant differential form of top degree $n$  extends holomorphically to a smooth  
compactification of $\calB/G$ if and only if 
the associated modular form is cuspidal (compare [Fr], Satz III.2.6).
\neupara{A special group}%
We equip $V=\cz^4$ with the hermitian form
$$\spitz{a,b}=\bar a_1b_2+\bar a_2b_1-\bar a_3b_3-\bar a_4b_4.$$
Let $\calE$ be the ring of Eisenstein numbers. We consider the lattice
$M=\calE^4$. We use the notation
$$G_3=\U(M).$$
We also consider the principal congruence subgroups
$$G_3[a]=\kernel\bigl(G_3\lo \GL(4,\calE/a)\bigr)\qquad (a\in\calE,\ a\ne 0).$$
In the paper [FSM] the ring of modular forms of  $G_3[3]$ has been determined.
We recall some of the main results. There are 15 forrms
$B_1,\dots,B_{15}$ of weight 1 and 10 cusp forms $C_1,\dots,C_{10}$ of weight
10. These 25 forms generate the ring of modular forms $A(G_3[3])$. The relations
between the generators have been determined:
$$\dim[G_3[3],k]=\cases{
0&  for $k<0$,\cr
1& for $k=0$,\cr
15& for $k=1$,\cr
130& for $k=2$,\cr
750& for $k=3$,\cr
3115&for $k=4$,\cr
-1377+(8019/2)k-2187tk^2+(729/2)k^3&for $k>4$.\cr}$$
and
$$\dim[G_3[3],k]-\dim[G_3[3],k]_0=\cases{
15& for $k=1$,\cr
120& for $k=2$,\cr
405& for $k=3$,\cr
765& for $k=4$,\cr
810& for $k>4$.}$$
Basic elements of $G_3$ are reflections along vectors $b$ of
norm $\spitz{b,b}=-1$. They are defined by
$$a\loma a-(1-\eta){\spitz{b,a}\over\spitz{b,b}}b$$
where $\eta$ is a 6th root of unity (hence a power of $-\zeta)$.
They transform $b$ to $\eta b$ and act as identity on the orthogonal
complement of $b$. Their order equals the order of $\eta$. We call them
biflections, triflections or hexflections corresponding to the order of
$\eta$.
\smallskip
By a short mirror in $\tilde\calB$ we understand
the fixed point set of a triflection. The short mirrors decompose under $G_3[\sqrt{-3}]$
into 15 orbits.
\medni
{
\halign{\qquad$#$\quad&$#$\quad&$#$\quad&$#$\quad&$#$\quad\cr
1&2&3&4&5\cr
(0,0,1,0)&(0,0,0,1)&(1,0,1,0)&(1,0,-1,0)&(1,0,0,1)\cr
\noalign{\vskip2mm}
6&7&8&9&10\cr
(1,0,0,-1)&(0,1,1,0)&(0,1,-1,0)&(0,1,0,1)&(0,1,0,-1)\cr
\noalign{\vskip2mm}
11&12&13&14&15\cr
(\zeta,-1,1,1)&(\zeta,-1,1,-1)&(\zeta,-1,-1,1)&(\zeta,-1,-1,-1)&(\zeta,1,0,0)\cr}
\medni
We recall the basic properties of the forms $B_1,\dots,B_{15}$. Their zero divisor
(considered mod $G_3[3]$) consist of three short mirrors which are described in the
following list:
\smallni
\vbox{\rm\noindent
(1,2,15),
(2,4,8),
(2,3,7),
(1,6,10),
(1,5,9),
(12,13,15),
(11,14,15),
(4,6,11),
(4,5,12),
(8,10,14),
(8,9,13),
(3,6,13),
(3,5,14),
(7,10,12),
(7,9,11).}
\smallni
We also recall the definition of the cusp forms $C_i$.
$$\eqalign{
C_1&=B_2B_4B_{15}/B_8,\cr
C_2&=B_2B_{13}B_{15}/B_3,\cr
C_3&=B_3B_6B_{10}/B_{14},\cr
C_4&=B_3B_5B_8/B_{15},\cr
C_5&=B_8B_{13}B_{14}/B_9,\cr
C_6&=B_5B_7B_{14}/B_{15},\cr
C_7&=B_2B_6B_{15}/B_{11},\cr
C_8&=B_1B_8B_{11}/B_2,\cr
C_9&=B_6B_{13}B_{15}/B_7,\cr
C_{10}&=B_2B_4B_6/B_1.\cr}$$
Each of them vanishes along 6 short mirrors (counted mod $G_3[3]$).
For example $C_1^2$ vanishes along the mirrors with the digits
$1,2,7,8,9,10$.
We denote by $G'$ the group which is generated by $G_3[3]$ and the 6 triflections
corresponding to the above 6 short mirrors. 
\smallskip
From \InvDif\ follows that $C_1^2$ defines a $G'$-invariant
differential form. Moreover its zeros (of order two) disappear in the quotient 
$\calB/G'$ due to the ramification. Hence
there is a chance that a suitable desingularization gives a  Calabi--Yau manifold.
We will show now that this is actually the case.
\proclaim
{Lemma}
{The index of the group $G_3[3]$ in $G'$ is $486$. It contains the negative of the identity,
The covering degree of $X(G_3[3])\to X(G')$ is $243$.}
IndGr%
\finishproclaim
{\it Proof.\/} The proof can be given by means of a direct computation with the help of
a computer.\qed
\smallskip
The forms $B_6$, $B_7$, $B_8$, $B_9$, $B_{12}$, $B_{13}$
do not vanish along one of the six short mirrors. Hence they have trivial multiplier
with respect to $G'$.
\proclaim
{Theorem}
{The six modular forms above generate the ring $A(G')$. Defining relations are
$$B_7B_9B_{12}=B_6B_8B_{13},\quad B_6^3+B_7^3-B_8^3+B_9^3-B_{12}^3-B_{13}^3=0.$$
The associated variety $\calX$ has  $108$
singularities which are nodes.}
LemS%
\finishproclaim
{\it Proof.\/} From the description of the zeros one can derive that the 6 forms have no joint zero.
By a result of Hilbert, the ring $A(G')$ is integral over the subring
$\cz[B_6,B_7,B_8,B_9,B_{12},B_{13}]$. The relations described in the theorem generate a prime
ideal. A dimension argument shows that they are defining relations between the 6 forms.
The Hilbert polynomial of the subring can be computed. The result is $(27t^2+9t^3)/3!$.
Its highest coefficient is $3/2$. We also know that the highest coefficient of the Hilbert polynomial
of $A(G')$ is $729/2$. The quotient of the two highest coefficients, $243$, equals the covering
degree $[G':G_3[3]]$ of $X(G')\to X(G_3[3])$. So we see that $A(G')$ and the subring
$\cz[B_6,B_7,B_8,B_9,B_{12},B_{13}]$ have the same field of fractions. So $A(G')$ is the normalization
of this subring. One can check that the relations in the theorem define a normal ring.
In fact it is a complete intersection  and regular in codimension 2.
This finishes the proof.\qed
\smallskip
We introduce variables $X_0,\dots,X_5$ and consider the homomorphism
$$\cz[X_0,\dots,X_5]\lo A(G')$$
which is defined by
$$X_0\loma -B_6,\ X_1\loma B_8,\  X_2\loma B_{13},\
X_3\loma B_7,\ X_4\loma B_9,\ X_5\loma -B_{12}.$$
Then the defining relations get
\rahmen{X_0X_1X_2=X_3X_4X_5,\quad X_0^3+X_1^3+X_2^3=X_3^3+X_4^3+X_5^3.}
This is a very well-known variety. We  learned from  [Me], Chap.~5, Sect.~6, 
that it has 108 singularities which are all nodes. They correspond
to certain cusps in the ball-model. Moreover there is a projective small resolution which is a rigid
Calabi-Yau manifold ($h^{12}=0$) with Euler number $72$. The Picard number is $h^{11}=36$.
\proclaim
{Theorem}
{The variety $\calX$ defined in $P^5\cz$ as intersection of two cubics
$$X_0X_1X_2=X_3X_4X_5,\quad X_0^3+X_1^3+X_2^3=X_3^3+X_4^3+X_5^3$$
can be identified with (compactified) ball quotient $X(G')$ with respect to the group $G'$.
This gives an example of a three dimensional ball-quotient which admits
a Calabi--Yau model.}
MTH%
\finishproclaim
There is an obvious group $H$ of automorphisms of order $5832$ acting on this variety.
First we can consider permutations of the variables $X_0,X_1,X_2$ and $X_3,X_4,X_5$
separately and by interchanging these two blocs. This gives a group of order $72$.
Then we consider transformations $X_\nu\mapsto \zeta_\nu X_\nu$
where $\zeta_\nu$ are third roots of unity with the property
$$\zeta_0\zeta_1\zeta_2=\zeta_3\zeta_4\zeta_4.$$
Taking into account that $(\zeta,\dots,\zeta)$ acts as identity, this gives us group
of order $81$ acting on the projective variety. Both types generate a group of order
$72\cdot 81=5832$ as has been stated. In fact, these
automorphisms are modular.
\proclaim
{Remark}
{Let $N'$ be the normalizer of $G'$ in the full modular group.
There is a natural surjective homomorphism
$$N'\lo H$$
which is compatible with the actions of $N'$ on $X(G')$ and of $H$ on $\calX$.
The kernel of this homomorphism
is generated by $G'$ and by the transformation
``multiplication by $\zeta$''.}
NormStrich%
\finishproclaim
{\it Proof.\/} Since we know the action of the full modular group on the generators, the statement can
by checked by means of a direct computation. We did it with the help of a computer.\qed
\smallskip
We want to determine the subgroup that leaves the Calabi--Yau form invariant.
The Calabi--Yau form (given in the modular picture by the modular form $C_1^2$) can be described
in the model $\calX$ as follows.
Consider the differential form
$$\omega={\displaystyle\sum_{i=0}^5 z_i\>dz_0\wedge\dots\wedge\widehat{dz_i}\wedge\dots\wedge dz_5\over
(z_0z_1z_2-z_3z_4z_5)(z_0^3+z_1^3+z_2^3-z_3^3-z_4^3-z_5^3)}$$
This is a meromorphic differential form on $P^5\cz$ with poles on the union of the two cubics.
Taking the Poincar\'e residue we get a  the Calabi--Yau form on $\calX$.
The group $H$ extends to $P^5\cz$ in an obvious way. It acts on $\omega$ with a character $\chi$.
This character can be computed easily:
\proclaim
{Remark}
{The character $\chi$ of $H$ that describes the action of $H$ on the Calabi--Yau form is
defined by the following properties:
\vskip1mm
\item{\rm 1)} For a permutation of $X_0, X_1,X_2$ and similarly of $X_3,X_4,X_5$ it is the sign
of the permutation.
\item{\rm 2)} For the permutation $X_0\leftrightarrow X_3$,
$X_1\leftrightarrow X_4$, $X_2\leftrightarrow X_5$ it is $-1$.
\item{\rm 3)} For a transformation
$X_\nu\mapsto \zeta_\nu X_\nu$
where $\zeta_\nu$ are third roots of unity with the property
$\zeta_0\zeta_1\zeta_2=\zeta_3\zeta_4\zeta_5$ the value of $\chi$ is
$\zeta_0\zeta_1\zeta_2$. The kernel of $\chi$ is a subgroup of index $6$ in $H$
which has the order $972$.\vskip0pt}
AutH%
\finishproclaim
It might be worthwhile to study subgroups of the kernel of $\chi$ 
such that the corresponding quotient of $\calX$
admits a Calabi--Yau model.
 \vskip1cm
\noindent
{\paragratit References}%
 \vskip1cm
\item{[CFS]} Cynk, S., Freitag, E., Salvati Manni, R.:
{\it The geometry and arithmetic of a Calabi-Yau Siegel threefold,\/}
Int.~Jour.~Math.\  {\bf 29}, 1561--1583  (2011)  (arXiv: 1004.2997)
\medni
\item{[Fr]} Freitag, E.: {\it Siegelsche Modulfunktionen,\/} Grundlehren
der mathematischen Wissenschaften, Bd.\ {\bf 254}. Berlin-Heidelberg-New
York: Springer (1983)
\medskip
\item{[FS1]} Freitag, E., Salvati Manni, R.:
{\it Some Siegel threefolds with Calabi-Yau model,\/}
Ann. Scuola Norm. Sup. Pisa Cl.\ Sci.\ (5). Vol IX, 833-850 (2010) (arXiv: 0905.4150)
\medni
\item{[FS2]} Freitag, E., Salvati Manni, R.:
{\it Some Siegel threefolds with Calabi-Yau model II,\/}
To appear on Kyungpook Math. Journal] (arXiv: 1001.0324)
\medni
\item{[FS3]} Freitag, E., Salvati Manni, R.:
{\it On Siegel three folds with a projective Calabi--Yau model,\/}
 Communications in Number Theory and Physics. Volume 5, Number 3	p.713-750 ( 2011)
(arXiv: 1103.2040 )
\medni
\item{[FS4]} Freitag, E., Salvati Manni, R.:
{\it A three dimensional ball quotient,\/}
preprint (2012)(arXiv: 1201.0131)
\medni
\item{[Me]} Meyer, C.: {\it A dictionary of modular threefolds,\/}
thesis, University of Mainz, Fachbereich Mathematik und Informatik (2005)
\bye